# The importance of finding the upper bounds for prime gaps in order to solve the twin primes conjecture and the Goldbach's conjecture

Andrea Berdondini

ABSTRACT. In this article we present a point of view that highlights the importance of finding the upper bounds for prime gaps, in order to solve the twin primes conjecture and the Goldbach's conjecture. For this purpose, we present a procedure for the determination of the upper bounds for prime gaps different from the most famous and known approaches. The proposed method analyzes the distribution of prime numbers using the set of relative numbers $\mathbb{Z}$. Using negative numbers too, it becomes intuitive to understand that that the arrangement of 2P+1 consecutive numbers that goes -P to P, is the only arrangement that minimizes the distance between two powers having the same absolute value of the base D, with $|D| \leq P$. This arrangement is considered important because by increasing the number of powers of the prime numbers within a range of consecutive numbers, it is presumed to decrease the overlap between the prime numbers considered. Consequently, by reducing these overlaps, we suppose to obtain an arrangement, in which the prime numbers less than and equal to P and their multiples occupy the greatest possible number of positions within a range of 2P+1 consecutive numbers. If this result could be demonstrated, would imply not only the resolution of the Legendre's conjecture, but also a step forward in the resolution of the twin primes conjecture and the Goldbach's conjecture.

**Introduction**

In this article, we will see the importance of finding the upper bounds for prime gaps or demonstrating the Legendre's conjecture [1] in order to solve the twin primes conjecture [2] and the Goldbach's conjecture [3]. The first step we must do is to reformulate the conjecture on the twin primes and the Goldbach's conjecture, so that they represent a problem very similar to that of finding the upper bound about the gap between two successive prime numbers. The conjecture of twin primes can be reformulated as follows: there are infinite twin prime numbers if for each number $P_{n1} = 6n_1 \pm 1$ ($n_1 \in \mathbb{N}, n_1 > 0$), the maximum numbers of consecutive natural numbers that can be written with this equation:

$$P_{n2}n_p \pm n_2, for\ n_p \in \mathbb{N}, n_p \geq n_2\ and\ P_{n2} = 6n_2 \pm 1, for\ n_2 \in \mathbb{N}, 0 < n_2 \leq n_1$$

is less than $P_{n1} + 2n_1 + 3$. Instead, Goldbach's conjecture can be reformulated as follows: each even number greater than two can be written as the sum of two primes if for each even number $P_a$ greater than 2 the maximum number of natural numbers greater than or equal to 2 and less or equal to $P_a - 2$, which can be written with these two equations:

$$Pn, for\ n \in \mathbb{N}, n > 1\ and\ 2 \leq P\ prime\ number < \sqrt{P_a}$$
$$P_a - Pn, for\ n \in \mathbb{N}, n > 1\ and\ 3 \leq P\ prime\ number < \sqrt{P_a}$$

is less than or equal to $Pa - 4$. In this way, the twin prime conjecture and the Goldbach's conjecture are very similar to the Legendre's conjecture. In practice all three conjectures are related to each other respect to the following problem: given a prime number P, what is the maximum number of consecutive positions that can occupy the prime numbers less than and equal to P and their multiples. This problem is fundamental in defining the upper bounds for prime gaps. Being able to solve this problem would most likely imply the resolution of the Legendre conjecture and a significant progress in solving the conjecture on the twin primes and the Goldbach's conjecture. Indeed, as can be seen, from the equations reported, the main difference between the Legendre's conjecture and that on the twin primes and that of Goldbach is that in the last two conjecture the prime numbers are translated by a constant.

Assuming that the arrangement of the prime numbers and their multiples, in which the greatest number of consecutive positions are occupied, is the arrangement where the prime numbers considered overlap each other as little as possible; we will analyze the arrangement that minimizes the distance between two powers having the same absolute value of the base D, with $|D| \leq P$. This type of analysis is done using the set of relative numbers $\mathbb{Z}$. We can use the set of relative numbers because we exploit the fact that, given prime number P, the prime numbers less than and equal to P create a pattern, in which all the possible arrangements of the considered prime numbers are present, which is repeated with a frequency $F = 2 \cdot 3 \cdot \ldots \cdot P$ obtained by multiplying P by the prime numbers less than P. Therefore given a prime number P, the frequency F will never be infinite, so we can develop a modular arithmetic of modulus F, in which the first terms are consecutive to the last terms. In practice the first 10 terms of this pattern go from 1 to 10, instead the last 10 terms go from -9 to 0. So the number zero represents F the last term of this pattern, in which all the prime numbers considered overlap. In this way we can pass from the set of natural numbers $\mathbb{N}$ to the set of relative integers $\mathbb{Z}$. Consequently, the minimum distance between two powers, having the same absolute value of the base D, is not $D - D^2$ but 2D (the distance between –D and D). So the arrangement of 2P+1 consecutive numbers in which two powers, having the same absolute value of the base D with $|D| \leq P$, are at the minimum distance is the one that goes from –P to P. The sequence going from –P to P is particularly interesting because it also contains the -1 and 1, two numbers that are not multiples of any prime number. So be able to prove that the arrangement in which there is the minimum distance between two powers, having the same absolute value of the base D with $|D| \leq P$, is also the arrangement where the prime numbers, less than and equal to P and their multiples, occupy the maximum number of positions on an interval containing 2P+1 consecutive numbers, would imply not only the resolution of the Legendre's conjecture, but also a step forward in the resolution of the twin primes conjecture and the Goldbach's conjecture.

**Reformulation of the twin primes conjecture**

We start by ordering the odd numbers using the arithmetic progressions $6n - 1$ e $6n + 1$ with $n \in \mathbb{N}$. These two arithmetic progression generate all pairs of twin prime numbers outside the pair formed by 3 and 5.

$$6n - 1, for\ n\ \in \mathbb{N}, n > 0 \tag{1}$$
$$6n + 1, for\ n\ \in \mathbb{N}, n > 0 \tag{2}$$

|       | $6n - 1$ | $6n + 1$ |
|-------|----------|----------|
| $n = 1$ | 5        | 7        |
| $n = 2$ | 11       | 13       |
| $n = \cdots$ | ….  | ….       |

If we want to remove all the composite numbers within the sequences (1) and (2), we must calculate the following three products:

$$(6n_1 - 1)(6n_2 + 1) \tag{3}$$
$$(6n_1 - 1)(6n_2 - 1) \tag{4}$$
$$(6n_1 + 1)(6n_2 + 1) \tag{5}$$

Let us start with the product (3) by setting $n_2 = n_1$.

$$(6n_1 - 1)(6n_1 + 1) = 6(6n_1^2) - 1$$

We note that from the product between two numbers belonging one to the sequence (1) and the other to the sequence (2) we obtain a number that falls into the sequence (1). Now we analyze the term $6n_1^2$, which determines the value of $n$ in the sequence (1).

$$6n_1^2 = (6n_1 - 1 + 1)n_1$$

Defining $P_{n1} = 6n_1 - 1$ we have:

$$(P_{n1} + 1)n_1 = P_{n1}n_1 + n_1$$

Being the common difference of the arithmetic progressions (1) and (2) obtained by the product between the prime numbers 2 and 3, this value will never be divisible by the other prime numbers greater than 3, therefore all the other products are at a distance $6P_{n1}$.

At this point we can deduce the formula that defines the values of $n$ in the sequence (1), in which there is a composite number obtained by multiplying $P_{n1}$ by a number $P_{n2} > P_{n1}$ generated by the sequence different from that which generates $P_{n1}$.

$$n = P_{n1}n_p + n_1, \text{ for } n_p \in \mathbb{N}, n_p \geq n_1 \text{ and } P_{n1} = 6n_1 \pm 1, \text{ for } n_1 \in \mathbb{N}, n_1 > 0 \tag{6}$$

Now let us consider the product (4) by setting $n_2 = n_1$.

$$(6n_1 - 1)(6n_1 - 1) = 1 - 12n_1 + 36n_1^2 = 6(-2n_1 + 6n_1^2) + 1$$

The value generated by this product falls in the sequence (2), therefore we analyze the term $-2n_1 + 6n_1^2$, which determines the value of $n$ in the sequence (2).

$$-2n_1 + 6n_1^2 = n_1(-2 + 6n_1)$$
$$n_1(-2 + 6n_1) = n_1(-1 + 6n_1 - 1)$$

Defining $P_{n1} = 6n_1 - 1$ we have:

$$n_1(-1 + 6n_1 - 1) = n_1(P_{n1} - 1)$$
$$n_1(P_{n1} - 1) = P_{n1}n_1 - n_1$$

We have shown that every square of every number generated by the sequence (1) is on the sequence $6n + 1$ where $n = P_{n1}n_1 - n_1$.

Knowing that all the other products are at a distance $6P_{n1}$, we can deduce the formula that defines the values of $n$ in the sequence (2), in which there is a composite number obtained by multiplying $P_{n1} = 6n_1 - 1$ by another number $P_{n2} = 6n_2 - 1$ with $P_{n2} \geq P_{n1}$.

$$n = P_{n1}n_p - n_1, for\ n_p \in \mathbb{N},\ n_p \geq n_1\ and\ P_{n1} = 6n_1 - 1, for\ n_1 \in \mathbb{N},\ n_1 > 0 \qquad (7)$$

Now let us consider the product (5) by setting $n_2 = n_1$.

$$(6n_1 + 1)(6n_1 + 1) = 1 + 12n_1 + 36n_1^2 = 6(2n_1 + 6n_1^2) + 1$$

The value generated by this product falls in the sequence (2), therefore we analyze the term $2n_1 + 6n_1^2$, which determines the value of $n$ in the sequence (2).

$$2n_1 + 6n_1^2 = n_1(2 + 6n_1)$$
$$n_1(2 + 6n_1) = n_1(1 + 6n_1 + 1)$$

Defining $P_{n1} = 6n_1 + 1$ we have:

$$n_1(1 + 6n_1 + 1) = n_1(P_{n1} + 1)$$
$$n_1(P_{n1} + 1) = P_{n1}n_1 + n_1$$

We have shown that every square of every number generated by the sequence (2) is on the sequence $6n + 1$ where $n = P_{n1}n_1 + n_1$.

Knowing that all the other products are at a distance $6P_{n1}$, we can deduce the formula that defines the values of $n$ in the sequence (2), in which there is a composite number obtained by multiplying $P_{n1} = 6n_1 + 1$ by another number $P_{n2} = 6n_2 + 1$ with $P_{n2} \geq P_{n1}$.

$$n = P_{n1}n_p + n_1, for\ n_p \in \mathbb{N},\ n_p \geq n_1\ and\ P_{n1} = 6n_1 + 1,\ for\ n_1 \in \mathbb{N}, n_1 > 0 \qquad (8)$$

<u>Taking into consideration each prime number $P_{n1}$ greater than 3 we can determine, using the formulas obtained, the position of each composite number present in the sequences (1) and (2).</u>

At this point we have found the formulas that determine the values of $n$, in which there will be at least an odd non-prime number generated by one of the two arithmetic progressions $6n - 1$ and $6n + 1$. So the values of $n$ in which there are no pairs of twin prime numbers are described by the formulas (6) (7) and (8), which can be grouped in the following formula:

$$n = P_{n1}n_p \pm n_1, for\ n_p \in \mathbb{N},\ n_p \geq n_1\ and\ P_{n1} = 6n_1 \pm 1, for\ n_1 \in \mathbb{N}, n_1 > 0 \qquad (9)$$

<u>So all the values of the sequence of natural numbers that cannot be written with the previous formula, will indicate a value of $n$ where the odd numbers generated by arithmetic progressions $6n - 1$ and $6n + 1$ are prime and create a pair of twin prime numbers.</u>

Analyzing the formula (9) we notice that every time we consider a greater number of the previous ones, the first value generated is $P_{n1}n_1 + n_1$ or $P_{n1}n_1 - n_1$. Consequently, since the arithmetic progressions (1) and (2) consist of infinite values, there will be infinite intervals in which only multiples of the numbers from 5 to $P_{n1}$ will be present. At this point if we could prove that there is always at least one number $n$, that we cannot write with equation (9), within these intervals, which we know to be infinite, we can prove that there are infinite twin prime numbers.

In sequences (1) and (2) the minimum distance between two prime numbers can be 2 if we consider the pair $6n - 1$ and $6n + 1$ or 4 if we consider the pair $6n + 1$ and $6(n + 1) - 1$. Taking into account the second pair, so as to have the longest interval, we can calculate the length of this interval by making the difference between the first two numbers that are generated by equation (9), when considering the numbers $P_{n1} = 6n_1 + 1$ and $P_{n2} = 6(n_1 + 1) - 1$.

$$P_{n2}(n_1 + 1) - n_1 - 1 - P_{n1}n_1 - n_1$$

Knowing that $P_{n2} = P_{n1} + 4$ we have:

$$(P_{n1} + 4)(n_1 + 1) - n_1 - 1 - P_{n1}n_1 - n_1$$
$$P_{n1} + 2n_1 + 3$$

<u>Then given a value of $P_{n1} = 6n_1 + 1$, if it is possible to demonstrate that the maximum number of the positions occupied consecutively by the values generated by equation (9), considering the prime numbers P greater than 3 and less than and equal to $P_{n1}$, is less than $P_{n1} + 2n_1 + 3$, knowing that there are infinitely many intervals of length $P_{n1} + 2n_1 + 3$, we can prove the conjecture on twin primes numbers.</u>
<u>As can be seen, the problem thus formulated is very similar to the Legendre's conjecture; indeed it is a matter of finding the maximum length of the sequence of consecutive values generated by equation (9). This formula represents nothing more than a translation of the prime numbers and their multiples.</u>

**Reformulation of the Goldbach's conjecture**

Let us now consider the Goldbach's conjecture, also in this case we use a procedure similar to that used for the twin prime conjecture. In practice, given an even number $P_a$, we analyze the values belonging to $\mathbb{N}$ ranging from 2 to $P_a - 2$.

2
3
…
$P_a - 2$

We remove from this sequence all the numbers that can be written with these two equations:

$$Pn, for\ n \in \mathbb{N},\ n > 1\ and\ 2 \leq P\ prime\ number < \sqrt{P_a} \qquad (10)$$
$$P_a - Pn, for\ n \in \mathbb{N},\ n > 1\ and\ 3 \leq P\ prime\ number < \sqrt{P_a} \qquad (11)$$

With the formula (10) we remove all the numbers that are not prime. With the formula (11) we do the same, but with the only difference that in this case we start from $P_a$. In practice, formula (11) translates the values generated by formula (10). Consequently, a number N belonging to $\mathbb{N}$ in the range from 2 to $P_a - 2$, which cannot be written with equations (10) and (11) implies the existence of two prime numbers whose sum from $P_a$.

$$N = P_1$$
$$N = P_a - P_2$$
$$P_1 = P_a - P_2$$
$$P_a = P_1 + P_2$$

So if it is possible to proof that, for every even number $P_a$ greater than 2, the maximum number of values belonging to ℕ greater than or equal to 2 and less than or equal to $P_a - 2$, which can be written with equations (10) and (11) is less than or equal to $P_a - 4$, we can prove Goldbach's conjecture.

Consequently, we arrive at a situation analogous to the twin primes conjecture and the Legendre's conjecture. In this case the problem to be solved is to find the arrangement where the multiples of the prime numbers and their translations occupy the greatest number of positions consecutively within the range from 2 to $P_a - 2$.

**Analysis of the distribution of prime numbers using the set of relative numbers**

As anticipated in the introduction, the prime numbers less than and equal to P generate a pattern, in which are present all the possible arrangements of the considered prime numbers, which is repeated with frequency F = 2 * 3 * ....... P obtained by multiplying P by the prime numbers less than P. This pattern is fundamental because it also contains the arrangement, in which the prime numbers less than and equal to P and their multiples occupy the maximum number of consecutive positions.

Since the frequency F is not infinite, we can develop a modular arithmetic of modulus F, in which the first terms are consecutive to the last ones. The first 10 terms of this pattern go from 1 to 10, instead the last 10 terms go from -9 to 0. It is interesting to note that the number zero represents F the last term of this pattern, in which all the prime numbers considered overlap. Therefore we consider relevant to study the distribution of prime numbers using the set of relative integers ℤ. Using also negative numbers we can define the following sequence.

$$-P \ldots \ldots \ldots \ldots \ldots \ldots -1\ 0\ 1 \ldots \ldots \ldots \ldots \ldots \ldots P \tag{12}$$

In which it is intuitive to understand how this sequence minimizes the distance between two powers having the same absolute value of the base D, with $|D| \leq P$.

In this arrangement the minimum distance between two powers, having the same absolute value of the base D with $|D| \leq P$, is not $D - D^2$ but $2D$. Indeed $-D$ and $D$ are two powers that have the same absolute value of the base, therefore their distance is $2D$, the least possible. The study of the distribution of powers is very important, because we want to find the arrangement in which the numbers less than or equal to P overlap each other as little as possible.

So the next step is to try to demonstrate that the arrangement (12) is also the arrangement, in which the prime numbers, less than and equal to P, occupy the maximum number of positions in an interval that contains 2P+1 consecutive numbers. In order to solve this important problem we will present a procedure that we believe is very promising.

Let us start by changing the arrangement (12) considering only the odd numbers. We thus obtain the following arrangement of P+1 odd consecutive numbers.

$$-P \ldots \ldots \ldots \ldots \ldots \ldots -1\ \ 1 \ldots \ldots \ldots \ldots \ldots \ldots P \tag{13}$$

We define two groups of odd numbers: $D_{ma}$ and $D_m$.

$$P/2 < D_{ma} \leq P$$

$$1 < D_m < P/2$$

Now we only consider the odd numbers $D_{ma}$, these numbers can at most be present twice inside the arrangement (13), which we know contain P+1 odd consecutive numbers.

Taking into consideration only the odd numbers $D_{ma}$ we try to find the arrangement, in a range consisting of P+1 odd consecutive numbers, in which the greatest possible number of positions are occupied. The arrangement that solves this problem is the arrangement (13).

The reason is that this arrangement is the only arrangement, in which all the odd numbers considered occupy two positions. Indeed, the prime number P, in order to occupy two positions within a range consisting of P+1 odd consecutive numbers, must occupy the first and last positions. Consequently, the odd number equal to P-1, must occupy the second and penultimate positions. Continuing iteratively for the other odd numbers, it is shown that the arrangement (13) is the arrangement of P+1 odd consecutive numbers, in which the largest number of positions are occupied considering the odd numbers less than or equal to P and greater than P/2 .

Now we take into consideration the odd numbers $D_m$, in this case different arrangements can exist compared to (13), in which these numbers occupy an extra position. Therefore, we try to understand what happens when we translate an odd number $D_m$ so that it occupies an extra position. In this case, there will always be a position occupied in the range from $D_m \lfloor P/D_m \rfloor$ to P or in the range from $-D_m \lfloor P/D_m \rfloor$ to -P. The reason is that the arrangement that goes from $-D_m \lfloor P/D_m \rfloor$ to $D_m \lfloor P/D_m \rfloor$ is the arrangement of $D_m + 1$ odd consecutive numbers, where $D_m$ occupies the maximum number of positions. Therefore, the extra position occupied must be in the range from $D_m \lfloor P/D_m \rfloor$ to P or in the range from $-D_m \lfloor P/D_m \rfloor$ to -P. Since the value $D_m \lfloor P/D_m \rfloor$ always greater than P/2, this implies that the extra position occupied by an odd number $D_m$ overlaps with an odd number $D_{ma}$. At this point, in order to keep the gain of the extra position, we will have to move the odd number $D_{ma}$, however, as shown above, there is only one arrangement in which each odd number $D_{ma}$ occupies two positions. Consequently, moving the odd number $D_m$ implies that the new arrangement, of the odd numbers $D_{ma}$, occupies one position less than the case of the arrangement (13).

So a $D_m$ number in order to occupy an additional position must necessarily occupy at least one position occupied by a $D_{ma}$ number, consequently the length calculated by the ends not occupied in the arrangement (13) is reduced, therefore a $D_{ma}$ number will occupy one position less. The reason is that the $D_{ma}$ numbers cannot occupy two positions if the distance, between the unoccupied ends in a range of odd consecutive numbers, is less than $2D_{ma}$.

We report the following example: if the last three positions in the arrangement (13) are occupied by the translation of the numbers $D_m$, the numbers: P, P-1 and P-2 will never occupy two positions, so we will lose three positions. Consequently, the translation of the numbers $D_m$ has as final result an arrangement where an equal or lesser number of positions will be occupied with respect to the sequence (13).

The argument just made applies to every odd number $D_m$, therefore we can presume that there is no other arrangement, of P+1 odd consecutive numbers, in which the odd numbers less than or equal to P occupy one position more than the arrangement (13).

## Consequences of the proof of the Legendre's conjecture on the twin primes conjecture and on the Goldbach's conjecture

If the sequence (12) represents the sequence in which the prime numbers, less than or equal to P, and their multiples occupy the greatest number of positions in a range of 2P+1 successive numbers, it means that given a number N there is always a prime number P greater than N and less than $N + 2\sqrt{N} + 1$, therefore the Legendre's conjecture is true. Now let us see the implications of this result on the twin primes conjecture and on the Goldbach's conjecture.

We start by considering the conjecture on twin prime numbers. The equation (9) does nothing but translate prime numbers greater than 3, which we can write in this way $P_n = 6n \pm 1$ ($n_1 \in \mathbb{N}, n_1 > 0$), of a constant equal to $n$. At this point just take the arrangement (13) and, given a prime number $P_n$, apply the translations to all odd numbers greater than 3 and less than or equal to $P_n$. The new arrangement contains $P_n + 2n + 1$ consecutive values, the reason is that $P_n$ is translated by $n$, and therefore to be contained twice, it needs an interval that contains no less than $P_n + 2n + 1$ odd consecutive numbers.

In order to understand the procedure used, we can overlap the arrangement (13) with the arrangement obtained by translating the numbers $P_n = 6n \pm 1$ of $\pm n$ positions ($P_n, \pm n$). Taking as an example $P_n = 7$ we have:

$$-9 \quad -7 \quad -5 \quad -3 \quad -1 \quad 1 \quad 3 \quad 5 \quad 7 \quad 9$$
$$-(7,+1) \quad -(5,+1) \quad -(7,-1) \quad -(5,-1) \quad -1 \quad 1 \quad (5,-1) \quad (7,-1) \quad (5,+1) \quad (7,+1)$$

In practice we translate the numbers $P_n = 6n \pm 1$, along the set of relative odd integers, by a number of positions equal to $n$. So for example: in the case of 5 $n$ is equal to 1, so we must translate it by one position forward +$n$ and one position back -$n$. When it is moved forward, the 5 occupies the position occupied by the value 7 in the arrangement (13). On the other hand, when it is moved backwards, the 5 occupies the position occupied by the value 3 in the arrangement (13). So the number 5 will occupy two positions, corresponding to the two translations +$n$ and –$n$.

The translated arrangement which has been compared with the arrangement (13) is that which, given a prime number $P_n = 6n \pm 1$ and $F = 5 \cdot 7 \ldots . P_n$ (F is obtained by multiplying $P_n$ by the prime numbers greater than 3 and less than $P_n$), goes from (F-$P_n$) / 2-n to (F+$P_n$) / 2+n, hence an arrangement containing $P_n + 2n + 1$ consecutive numbers.

$$\frac{F-P_n}{2} - n \ldots\ldots\ldots\ldots\ldots\ldots . \frac{F-1}{2} \quad \frac{F+1}{2} \ldots\ldots\ldots\ldots\ldots . \frac{F+P_n}{2} + n \qquad (14)$$

In the case of P = 7 we have:

$$13 \quad 14 \quad 15 \quad 16 \quad 17 \quad 18 \quad 19 \quad 20 \quad 21 \quad 22$$

The central values of the shifted arrangement are (F-1)/2 and (F+1)/2 which are equivalent to -1 and 1 in the arrangement (13). These values will never be occupied by any value generated by equation (9). Consequently if the arrangement (14) is the arrangement where the values generated by the equation (9), considering the prime numbers from 5 to $P_n$, occupy the greatest number of positions on a range of $P_n + 2n + 1$ consecutive numbers, we proof the conjecture on twin prime numbers. Indeed we know, as demonstrated previously, that there are an infinite number of distinct intervals of $P_n + 2n + 3$ consecutive numbers, in which there are only the values generated by equation (9) when considering the prime numbers from 5 to $P_n$.

Now let us take Goldbach's conjecture into consideration, in this case we have a slightly more complex situation, because the value of the translations depends on the even number $P_a$. Given an even number $P_a$, equation (11) can be rewritten in this way:

$$Pn - \left(P\left(\left\lfloor\frac{P_a}{P}\right\rfloor + 1\right) - P_a\right) \text{ for } n \in \mathbb{N}, \quad 0 < n < \left\lfloor\frac{P_a}{P}\right\rfloor \text{ and } 3 \leq P \text{ prime numebr} < \sqrt{P_a} \quad (15)$$

At this point, in order to simplify the discussion, and to obtain a situation very similar to that just described with respect to the twin primes conjecture, we can remove the 3 and all its multiples and its relative translations calculated with equation (15), in which we also eliminates the value for $n = \lfloor P_a/P \rfloor$. In this way we remove 1/3 of all numbers between 2 and $P_a - 2$, obtaining a sequence of values spaced 3 from each other. From this range of values we must remove all multiples of the prime numbers greater than 3 and less than $\sqrt{P_a}$ and the relative translations. The values that need to be removed are generated by the following two equations:

$$Pn \text{ for } n \in \mathbb{N}, \quad n > 1 \text{ and } 3 < P \text{ prime numebr} < \sqrt{Pa}$$
$$Pn - \left(P\left(\left\lfloor\frac{P_a}{P}\right\rfloor + 1\right) - P_a\right) \text{ for } n \in \mathbb{N}, \quad 0 < n < \left\lfloor\frac{P_a}{P}\right\rfloor \text{ and } 3 < P \text{ prime number} < \sqrt{P_a}$$

In this way we obtain a situation very similar to that just described, where for each prime number $P_n = 6n \pm 1$, the translations $+n$ and $-n$ are replaced by 0 and $-(P(\lfloor P_a/P \rfloor + 1) - P_a)$. In order to obtain this analogy, very restrictive conditions have been imposed, in which we remove more values than those defined by equations (10) and (11).

This was done because the purpose of this discussion is not to give a mathematical proof of the conjecture, but only to bring attention to how easy it is to use this procedure also for problems other than the Legendre's conjecture. Consequently, if it were possible to proof that the arrangement (12) is the arrangement where the prime numbers, less than and equal to P, and their multiples occupy the maximum number of positions on an interval that contains 2P+1 consecutive numbers, it would be easy to apply this method on the twin primes conjecture and on the Goldbach's conjecture.

**Conclusion**

In this article, we have exposed a point of view that highlights the importance of finding the upper bounds for prime gaps and therefore solving the Legendre's conjecture, in order to solve the twin primes conjecture and the Goldbach's conjecture.

We have also analyzed a procedure for the determination of the upper bounds for prime gaps different from the more famous and known approaches [4], [5] and [6]. The proposed method analyzes the distribution of prime numbers using the set of relative numbers $\mathbb{Z}$. Using negative numbers, it becomes intuitive to understand that the arrangement (12) is the only arrangement, of 2P+1 consecutive numbers, which minimizes the distance between two powers having the same absolute value of the base D, with $|D| \leq P$.

The arrangement (12) is considered important because by increasing the number of powers of the prime numbers within a range of consecutive numbers, it is presumed to decrease the overlap between the prime numbers considered. Consequently, by reducing these overlaps, we suppose to obtain an arrangement, in which the prime numbers less than and equal to P and their multiples occupy the greatest possible number of positions within a range of 2P+1 consecutive numbers. This result, as explained in the previous chapters, is fundamental not only for solving the Legendre's conjecture but also for the twin primes conjecture and the Goldbach's conjecture.

*E-mail address*: andrea.berdondini@libero.it